\documentclass{article}
\usepackage[numbers]{natbib}
\usepackage{optml}
\usepackage{amsmath}

\usepackage{amsfonts}       
\usepackage{nicefrac}       
\usepackage{microtype}      
\usepackage{amsmath}
\usepackage{amssymb}
\usepackage{bm}
\usepackage{bbm}
\usepackage{esint}
\usepackage[pdftex]{graphicx}
\graphicspath{{results/}}
\usepackage{epstopdf}
\usepackage{subcaption}
\usepackage[ruled,vlined]{algorithm2e}
\DontPrintSemicolon
\usepackage[standard]{ntheorem}
\usepackage{cleveref}
\usepackage{mathtools}
\usepackage{wrapfig}

\newcommand{\calG}{\mathcal{G}}

\newcommand{\calB}{\mathcal{B}}
\newcommand{\calN}{\mathcal{N}}

\newcommand{\one}{\bm{1}}

\newcommand{\RR}{\mathbb{R}}

\newcommand{\prox}{\mathrm{prox}}

\newcommand{\bu}{\bm{u}}
\newcommand{\bx}{\bm{x}}
\newcommand{\by}{\bm{y}}

\newcommand{\bs}{\bm{s}}
\newcommand{\bb}{\bm{b}}

\newcommand{\bbeta}{\bm{\beta}}
\newcommand{\bt}{\bm{t}}
\DeclareMathOperator*{\argmax}{\mathrm{argmax}}
\DeclareMathOperator*{\argmin}{\mathrm{argmin}}

\newcommand{\dom}{\mathrm{dom}\,}

\newcommand{\bv}{\bm{v}}

\newcommand{\bA}{\bm{A}}

\newcommand{\APGncp}{APGnc\textsuperscript{+}}
\newcommand{\BCoAPGncp}{BCoAPGnc\textsuperscript{+}}

\newcommand{\sign}{\mathrm{sign}}

\newtheorem{assumption}{Assumption}
\newtheorem{notation}{Notation}

\title{Accelerated Block Coordinate Proximal Gradients with Applications in High Dimensional Statistics}
\author{\name Tsz Kit Lau \email{tklauag@connect.ust.hk}\\
  \addr{Department of Mathematics, Hong Kong University of Science and Technology, Hong Kong}\\
  \name Yuan Yao \email{yuany@ust.hk}\\
  \addr{Department of Mathematics, Hong Kong University of Science and Technology, Hong Kong}
}

\begin{document}
\maketitle

\begin{abstract}
  Nonconvex optimization problems arise in different research fields and arouse lots of attention in signal processing, statistics and machine learning. In this work, we explore the accelerated proximal gradient method and some of its variants which have been shown to converge under nonconvex context recently. We show that a novel variant proposed here, which exploits adaptive momentum and block coordinate update with specific update rules, further improves the performance of a broad class of nonconvex problems. In applications to sparse linear regression with regularizations like Lasso, grouped Lasso, capped $\ell_1$ and SCAD, the proposed scheme enjoys provable local linear convergence, with experimental justification. 
\end{abstract}

\section{Introduction}\label{sec:intro}
Many problems in machine learning are targeted to solve the following minimization problem 
\begin{equation}\label{min}
\min_{\bx\in\RR^n} F(\bx) \equiv f(\bx)+g(\bx), 
\end{equation}
where $f$ is differentiable, $g$ can possibly be nonsmooth. Convexity is not assumed for $f$ and $g$. If $F$ is convex, it is shown that \eqref{min} can be solved efficiently by the accelerated proximal gradient (APG) method (sometimes referred to as FISTA \cite{Beck2009}, as in \Cref{apg}. APG has a convergence rate of $O(1/k^2)$ which meets the theoretical lower bound of first-order gradient methods for minimizing smooth convex functions. 
 \begin{wrapfigure}{L}{0.5\textwidth}
     \begin{minipage}{0.5\textwidth}
 \begin{algorithm}[H]
 	\caption{APG}
 	\label{apg}
 	\SetKwInOut{Input}{Input}
 	\Input{$\by_1=\bx_1=\bx_0, t_1=1, t_0=0, \eta\le\nicefrac{1}{L}$. }
 	\For{$k=1, 2, \ldots$}{
 		$\by_k=\bx_k+\frac{t_{k-1}-1}{t_k}(\bx_k-\bx_{k-1})$,\;
 		$\bx_{k+1}=\prox_{\eta g}(\by_k-\eta\nabla f(\by_k))$,\;
 		$t_{k+1}=\frac{1+\sqrt{1+4t_k^2}}{2}$.
 	}
 \end{algorithm}
 \end{minipage}
   \end{wrapfigure}

For nonconvex version of \eqref{min}, APG was first introduced and analyzed by \citet{Li2015}, in which they propose monotone APG (mAPG) and nonmonotone APG (nmAPG) by exploiting the Kurdyka-\L{}ojasiewicz (KL) property. The main deficiency of these two algorithms is that they require two proximal steps in each iteration. In view of this, \citet{Yao2017} propose nonconvex inexact APG (niAPG) which is equivalent to APG for nonconvex problems (APGnc) in \cite{Li2017} if the proximal step is exact. This is an algorithm with comparable numerical performance to other state-of-the-art algorithms.

\citet{Li2017} analyze the convergence rates of mAPG and APGnc by exploiting the KL property. They further propose an \APGncp\,algorithm, which improves APGnc by introducing an adaptive momentum (\Cref{apgncp}, see Appendix). \APGncp\,has the same theoretical convergence rate as APGnc but has better numerical performance. However, the aforementioned APG-like algorithms do not leverage the special structure of the objective function $F$ (particularly the structure of $g$), which is what we want to explore in this paper.

In practice, many machine learning and statistics problems in the form of \eqref{min} have separable or block-separable regularizer $g$, so we can rewrite (\ref{min}) as 
\begin{equation}\label{blockmin}
\min_{\bx\in\RR^n} F(\bx) \equiv f(\bx)+\sum_{i=1}^{s} g_i(\bx_i), 
\end{equation}
where the variable $\bx=(\bx_1, \ldots, \bx_s)\in\RR^n$ has $s$ blocks, $s\ge 1$, and again the functions $g_i$, $i=1, \ldots, s$, can be nonconvex and nonsmooth. Block coordinate update is widely applied to solve convex and nonconvex problems in the form of \eqref{blockmin}. Since each iteration has low computational cost and small required memory, it is easy for parallel and distributed implementations and thus regarded as a more feasible method for large-scale problems. \citet{Xu2017} propose a block prox-linear (BPL) method (\Cref{bpl}, see Appendix) which can be viewed as a block coordinate version of APG. At iteration $k$, only one block $b_k\in\{1, \ldots, s\}$ is selected and updated. They establish the whole sequence convergence of BPL to a critical point, first by obtaining subsequence convergence followed by exploiting the KL property again. In their numerical tests, they mainly resort to randomly shuffling of the blocks and show that it leads to better numerical performance, as opposed to cyclic (Gauss-Seidel iteration scheme) updates \cite{Attouch2010,Attouch2013} and randomized block selection \cite{Lin2015}. 

\paragraph{Contribution.} Our main contribution is to bring adaptive momentum and block prox-linear method together in a new algorithm, for a further numerical speed-up sharing the same theoretical convergence guarantee under the KL property. Moreover, a new block update rule based on Gauss-Southwell Rule, is shown to beat randomized or cyclic updates as it selects the ``best'' block maximizing the magnitude of the step at each iteration. In the applications to high dimensional statistics, we particularly show that for sparse linear regressions with regularizations including (grouped) Lasso, capped $\ell_1$, and SCAD, the proposed algorithm has \emph{provable local linear convergence}. \S \ref{sec:alg} presents the algorithm and \S \ref{sec:app} discusses its applications in those sparse linear regressions together with empirical justification.

\section{The \BCoAPGncp\,Algorithm} \label{sec:alg}
We present our proposed algorithm in \Cref{bcapgncp}, Block-Coordinate APGnc with adaptive momentum (\BCoAPGncp), which takes advantage of several acceleration tools in \Cref{apgncp,bpl}. 

\begin{algorithm}[h]
	\caption{Block-Coordinate APGnc with adaptive momentum (\BCoAPGncp)}
	\label{bcapgncp}
	\SetKwInOut{Input}{Input}
	\Input{no. of blocks $s$, $\by^{(1)}=\bx^{(1)}=\bx^{(0)}=\bx^{(-1)}=\bv^{(0)}, \beta, t\in(0, 1), \bbeta^{(0)}=\beta\one_s, \bt=t\one_s$. }
	\For{$k=1, 2, \ldots$}{
		Pick $b_k\in\{1, \ldots, s\}$ in a deterministic (e.g., GS-$r$ rule) or random manner and set $\alpha_i^{(k)}$,\;
		$\begin{dcases*}
		\bx_i^{(k)} = \bx_i^{(k-1)} & if $i\ne b_k$,\\
		\bx_i^{(k)}\in\prox_{\alpha_i^{(k)}g_i}\left(\widehat{\bx}_i^{(k)} -\alpha_i^{(k)}\nabla_{\bx_i} f\left(\bx_{\ne i}^{(k-1)}, \widehat{\bx}_i^{(k)}  \right)  \right) & if $i=b_k$,\\
		\text{where }\widehat{\bx}_i^{(k)} :=\bx_i^{(k-1)}+\beta_i^{(k-1)}(\bx_i^{(k-1)}-\bx_i^{(k-2)}),
		\end{dcases*}$	\;
		$\begin{dcases*}
		\bv_i^{(k)} = \bv_i^{(k-1)} & if $i\ne b_k$,\\
		\bv_i^{(k)}=\bx_i^{(k)}+\beta_i^{(k-1)}(\bx_i^{(k)}-\bx_i^{(k-1)}) & if $i=b_k$,
		\end{dcases*}$\;
		\If{$F(\bx^{(k)})\le F(\bv^{(k)})$}{$\begin{dcases*}
				\beta_i^{(k)} =\beta_i^{(k-1)} & if $i\ne b_k$,\\
			\beta_i^{(k)}=\beta_i^{(k-1)}t_i &  if $i=b_k$,
		\end{dcases*}$}
		\ElseIf{$F(\bx^{(k)})\ge F(\bv^{(k)})$}{$\begin{dcases*}
		\beta_i^{(k)} =\beta_i^{(k-1)}	& if $i\ne b_k$,\\
			\beta_i^{(k)}=\min\left\lbrace \beta_i^{(k-1)}/t_i,1\right\rbrace &  if $i=b_k$. 
			\end{dcases*}$}
	}
\end{algorithm}
In each iteration of this algorithm, the updates of $\widehat{\bx}_i^{(k)}$ and $\bx_i^{(k)}$ follow that of \Cref{bpl}. The extrapolation step $\bv_i^{(k)}$ aims to further exploit the opportunity of acceleration by magnifying the momentum $\beta_i$ of the block $i$ when $\bv^{(k)}$ achieves an even lower objective value. If this does not hold, the momentum is diminished. This intuitive yet efficient step follows the main idea of \Cref{apgncp}. 
\paragraph{Gauss-Southwell rules.}
In \cite{Nutini2015}, three proximal-gradient Gauss-Southwell rules are presented, namely the GS-$s$, GS-$r$ and GS-$q$ rules (see Appendix for details). In particular, in \S \ref{sec:app}, we test with the GS-$r$ which would possibly speed up the convergence since it selects the block which maximizes the magnitude of the step at each iteration. 

We need certain assumptions of Problem \eqref{blockmin} in order to establish the convergence of our \Cref{bcapgncp}. 

\begin{assumption}\label{assumption}
	We have several assumptions on the functions $F, f$ and $g_i$: (i) $F$ is proper and bounded below in $\dom F$, $f$ is continuously differentiable, and $g_i$ is proper lower semicontinuous for every $i$. Problem \eqref{bpl} has a critical point $\bx^*$, i.e., $\bm{0}\in\partial F(\bx^*)$; (ii) Let $i=b_k$, then $\nabla_{\bx_i} f(\bx_{\ne i}^{(k-1)}, \bx_i)$ has Lipschitz constant $L_k$ with respect to $\bx_i$, which is bounded, finite and positive for all $k$; (iii) In \Cref{bcapgncp}\, every block is updated at least once within any $T$ iterations (see Appendix for definitions and notations). 
\end{assumption}
\begin{theorem}[Whole sequence convergence]\label{conv}
Suppose that \Cref{assumption} holds. Let $\{\bx^{(k)}\}_{k\ge 1}$ be generated from \Cref{bcapgncp}. Assume (i) $\{\bx^{(k)}\}_{k\ge 1}$ has a finite limit point $\bar{\bx}$; (ii) $F$ satisfies the KL property (\Cref{klprop}, see Appendix) around $\bar{\bx}$ with parameters $\rho$, $\eta$ and $\theta$; (iii) For each $i$, $\nabla_{\bx_i}f(\bx)$ is Lipschitz continuous within $B_{4\rho}(\bar{\bx})$ with respect to $\bx$. Then, we have $\lim_{k\to\infty}\bx^{(k)}=\bar{\bx}$. 
\begin{proof}
This theorem mainly follows from Theorem 2 of \cite{Xu2017}, since in its proof, there are no strict requirements on the momentum as long as it is between 0 and 1, and on the block update rule for each iteration. The step size $\alpha_k$ can be set to fulfill the assumption of this theorem which depends on the choice of Lipschitz constant $L_i$ of $\nabla_{\bx_i}f(\bx_{\ne i}^{(k-1)},\bx_i)$. 
\end{proof}	
\end{theorem}

\section{Applications in High Dimensional Statistics} \label{sec:app}
Many optimization problems in machine learning and statistics can be formulated in the form of \eqref{min}, for instance, sparse learning \cite{Bach2012}, regressions with nonconvex regularizers \cite{Fan2001,ZhangCH2010}, capped $\ell_1$-norm \cite{Zhang2010} and the log-sum-penalty \cite{Candes2008}. We consider the general class of regularized least squares or regression problems having the form 
\begin{equation}\label{reg_ls}
\min_{\bx\in\RR^n} \dfrac{1}{2n}\|\bA\bx-\bb\|^2+R(\bx), 
\end{equation}
where $\bA\in\RR^{m\times n}$, $\bb\in\RR^m$ and $R(\bx)$ is a separable or block separable regularizer. We test our algorithm with two convex problems and two nonconvex problems, mainly with sparse instances since GS rules can be efficiently calculated in such cases. These strategies can be applied to nonconvex problems since the calculation doe not depend on convexity \cite{Nutini2015}. All objective functions in this section are KL functions since all of them are semialgebraic \cite{Attouch2010,Xu2013}. 

\paragraph{$\ell_1$-regularized Underdetermined Sparse Least Squares.}
In this case, $R(\bx)=\lambda\|\bx\|_1$ which is a nonsmooth regularizer promoting sparsity. We generate $\bA$ and $\bb$ in the same way as that of \S9 in \cite{Nutini2015}, with $\lambda=1$, $m=1000$ and $n=5000$. In BPL and \BCoAPGncp, for illustrative purpose, we separate $\bx$ into $s$ blocks (set $s=5$) of equal size $N:=n/s$, and the step size at each iteration for a selected block $b_k$ is chosen to be $\alpha_i^{(k)}:=1/\|\bA_{:,b_k}^\top \bA_{:,b_k}\|_2$, where $\bA_{:,b_k}$ represents the columns of $\bA$ corresponds to the block $b_k$. We also take $\beta=t=0.9$ for \APGncp\,and \BCoAPGncp. 

\paragraph{$\ell_1/\ell_2$-regularized Sparse Least Squares.}
In this case, $R(\bx)=\lambda\sum_{g\in\calG} \|\bx_g\|_2$, where $\calG$ is a partition of $\{1, \ldots, n\}$ containing $\{1+kn/s, \ldots, (k+1)n/s\}$, where $k=0, \ldots, s-1$, for $s=5$ and $n=5000$. We use the same data set and parameters except $\beta=0.8$ and $t=0.2$. Only block coordinate methods are used since $R(\bx)$ is not completely separable. 

\paragraph{Capped $\ell_1$-regularized Sparse Least Squares.}
We consider the nonconvex capped $\ell_1$ penalty \cite{Zhang2010}, $R(\bx)=\lambda\sum_{i=1}^n\min\{|x_i|,\theta\}$, $\theta>0$. In this case we specify $s=10$, $\lambda=0.0001$, $\theta=0.1\lambda$, $\beta=0.8$ and $t=0.2$. 

\paragraph{Least Squares with SCAD penalty.}
We consider another nonconvex penalty term, the smoothly clipped absolute deviation (SCAD) penalty \cite{Fan2001}, $R(\bx)=\sum_{i=1}^n r_{\lambda, \gamma}(x_i)$, where $r_{\lambda, \gamma}(u)$ is defined in Appendix. Both $\bA$ and $\bb$ are sampled from $\calN(0,1)$ but they are standardized such that $\bb$ and each column of $\bA$ have zero mean, and each column of $\bA$ has unit variance. We take $m=1000$, $n=5000$, $s=10$, $\lambda=0.0001$ and $\gamma=3$. 

\begin{theorem}[Convergence rate]\label{conv_rate}
	Suppose $F$ is in the form of \eqref{reg_ls}, where $R(\bx)$ is chosen as the above four examples. Under the assumptions of \Cref{conv}, we have $\|\bx^{(k)}-\bar{\bx}\|\le C\alpha^k$, $\forall k$, for a certain $C>0$, $\alpha\in[0,1)$. Thus, $\{\bx^{(k)}\}_{k\ge 1}$ converges locally linearly to a stationary point of $F$.
\begin{proof}
According to Propositions 4.1, 4.2 and 5.2 of \cite{Li2017b}, all four above examples of $F$ are KL functions with an exponent $\theta=\nicefrac{1}{2}$ (for $\ell_1/\ell_2$-regularized least squares, we also need a mild assumption that $\inf_{\bx\in \RR^n} F(\bx)>\inf_{\bx\in \RR^n} \|\bA\bx-\bb\|^2/(2n)$). Then, the desired result follows immediately from Theorem 3 of \cite{Xu2017}. The convergence rate theorem for general $F$ is the same as this Theorem 3 for the same reason as in the proof of \Cref{conv}. 
\end{proof}	
\end{theorem}

We plot the value $F(\bx)-F(\bx^*)$ in each experiment, for the proposed algorithm and some existing APG-like algorithms mentioned in \S \ref{sec:intro}. For fair comparison, since in each iteration block coordinate methods only update one block, we consider $s$ block updates as one iteration in the plots. 
\begin{figure}[h!]
	\begin{subfigure}[h]{0.245\linewidth}
		\centering
		\includegraphics[height=3.8cm]{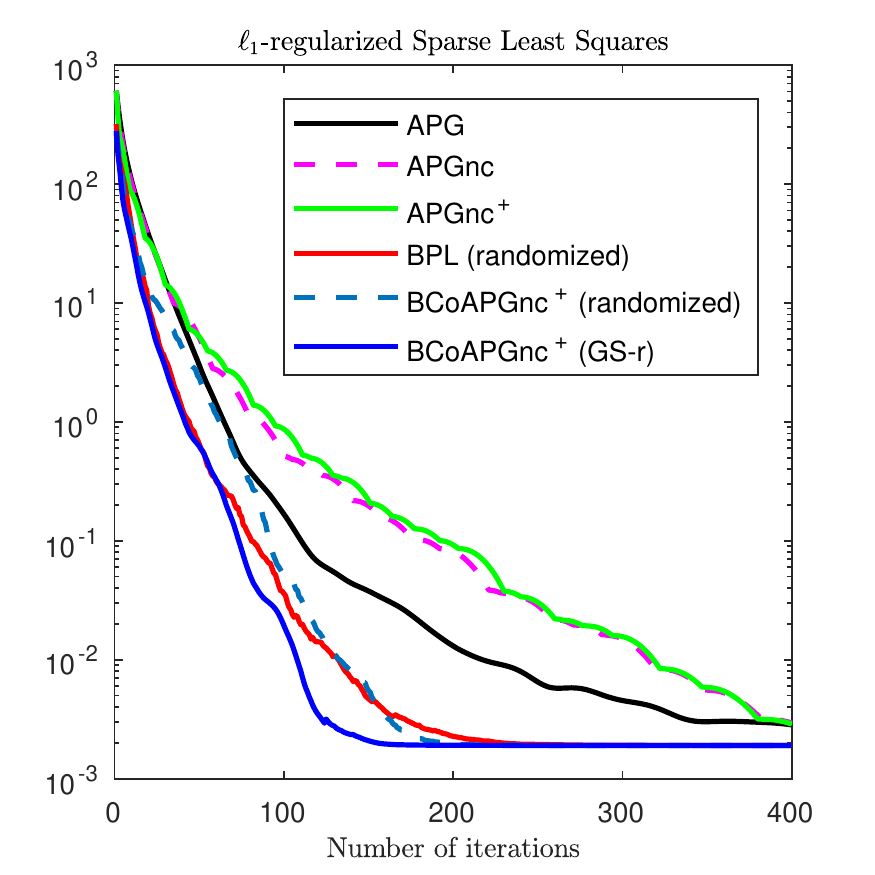}
		\caption{$\ell_1$ sparse LS}
		\label{fig1a}
	\end{subfigure}\hspace*{0.2mm}
	\begin{subfigure}[h]{0.245\linewidth}
		\centering
		\includegraphics[height=3.8cm]{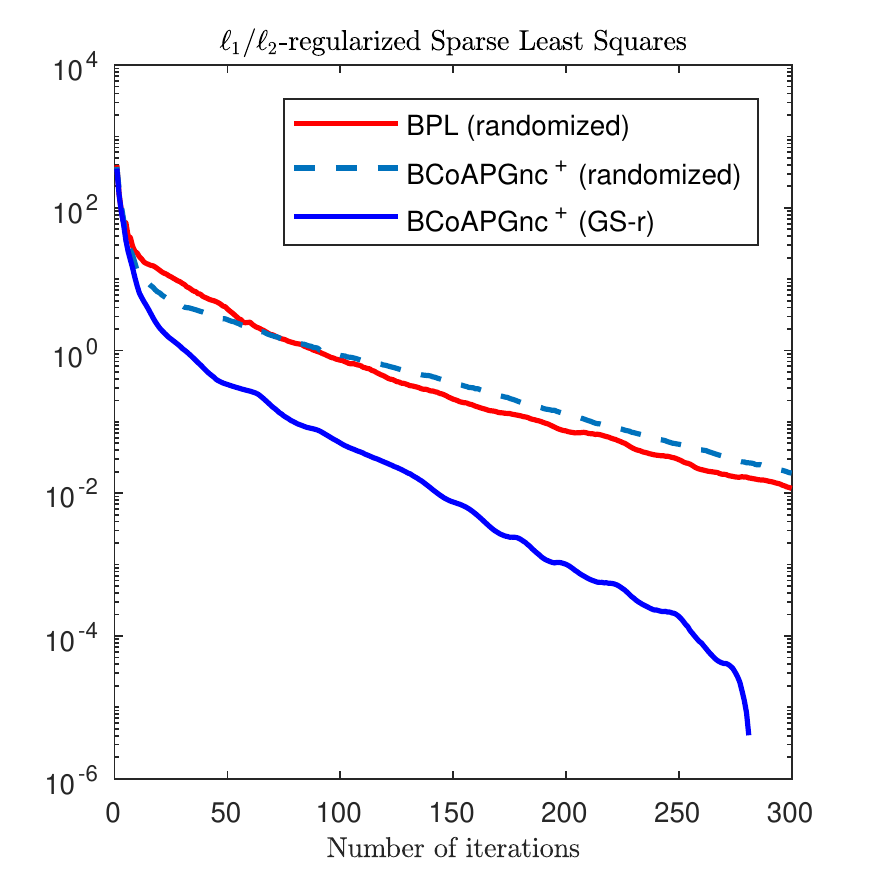}
		\caption{$\ell_1/\ell_2$ sparse LS}
		\label{fig1b}
	\end{subfigure}\hspace*{0.3mm}
	\begin{subfigure}[h]{0.245\linewidth}
		\centering
		\includegraphics[height=3.75cm]{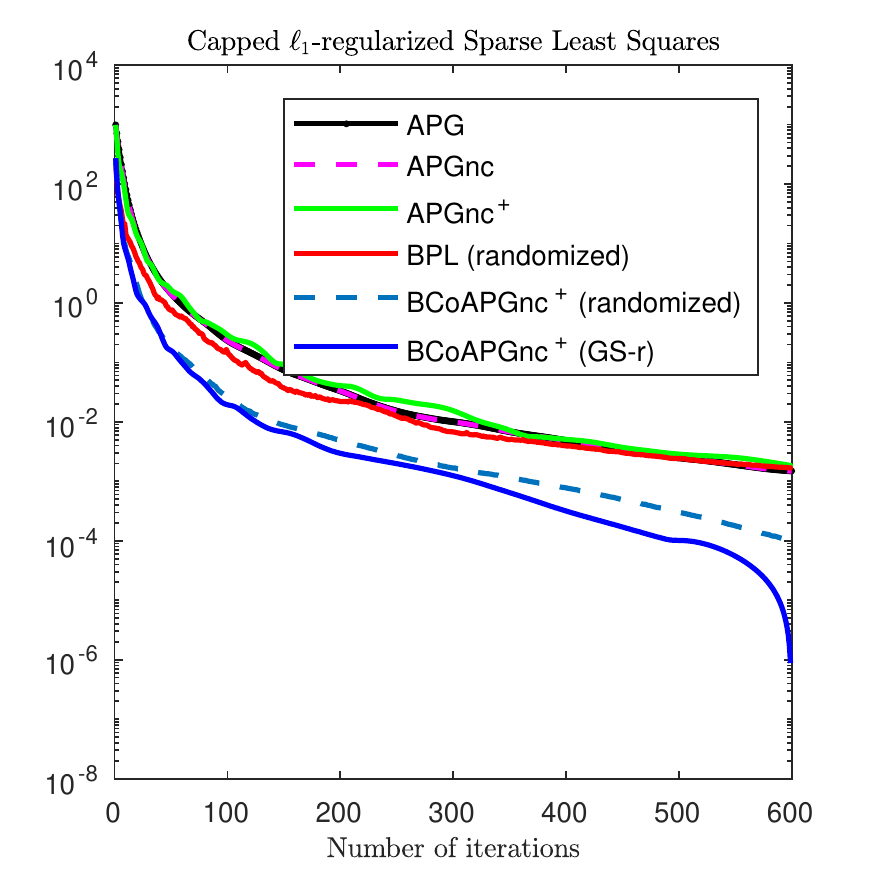}
		\caption{Capped $\ell_1$ sparse LS}
		\label{fig1c}
	\end{subfigure}
\begin{subfigure}[h]{0.245\linewidth}
	\centering
	\includegraphics[height=3.8cm]{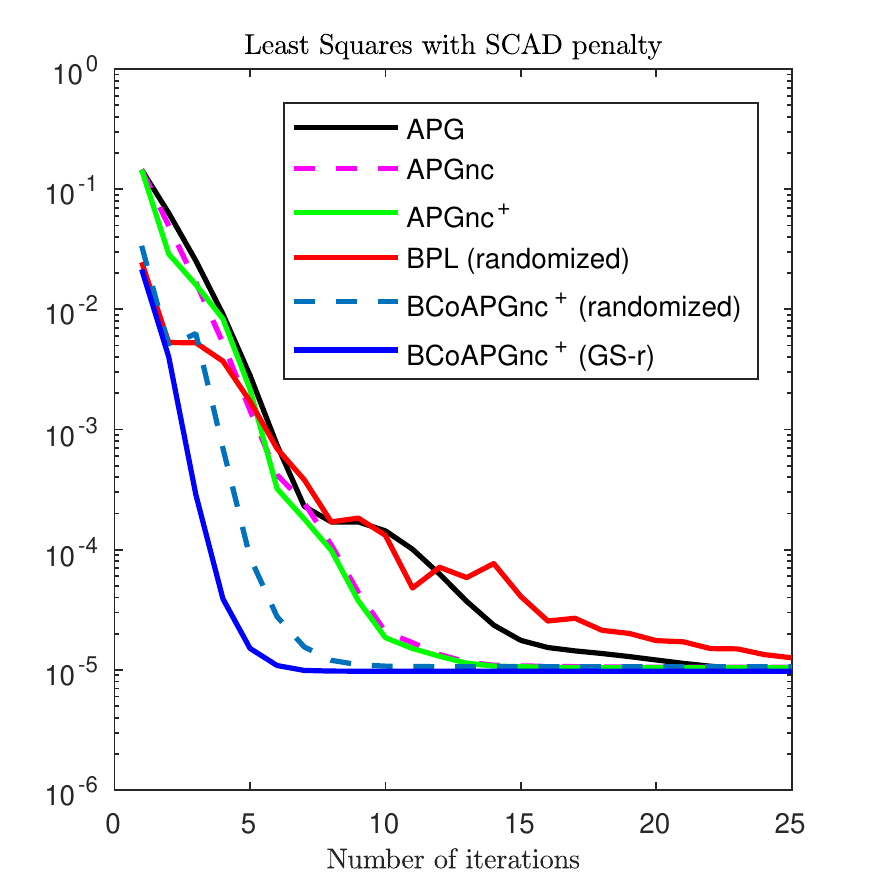}
	\caption{LS with SCAD}
	\label{fig1d}
\end{subfigure}
	\caption{Comparison of APG-like algorithms}
	\label{fig1}
\end{figure}
We observe in \Cref{fig1} that our proposed algorithm \BCoAPGncp\,(both randomized and GS-$r$ versions) provides the greatest initial acceleration, in both convex and nonconvex examples. It dominates most existing methods, especially during the first 20 iterations.  We also see that in general the \BCoAPGncp\,with GS-$r$ updates outperforms that with randomized updates, justifying the use of GS-$r$ rule for further acceleration. We see that \BCoAPGncp\,has superior performance in nonconvex problems revealed in \Cref{fig1c,fig1d}, where its original counterpart BPL does not give monotone objective value decline (in which its momentum is chosen according to that of APG). Overall, both versions of our proposed \BCoAPGncp\,speed up the convergence, compared with the current state-of-the-art \APGncp and BPL. 

\section{Conclusion and Future Work} \label{sec:con}
In this paper, we suggest a new algorithm which considers three main acceleration techniques, namely adaptive momentum, block coordinate update and GS-$r$ update rule. We also implement with adaptive step sizes in our experiments. We show that it shares the same convergence guarantee and convergence rate as BPL. It is noteworthy that all objective functions in sparse linear regression considered here, including (grouped) Lasso, capped $\ell_1$ and SCAD regularizations, are KL functions with an exponent $\nicefrac{1}{2}$, and thus have local linear convergence to their stationary points. Experiments show impressive results that our proposed method outperforms the current state-of-the-art. 

We focus on experiments with convex losses and (block) separable regularizers. Other applications with block separable regularizers but nonconvex losses such as matrix factorization and completion (e.g., in \cite{Xu2013}) deserve further treatment by applying the proposed algorithm. Further acceleration would be the use of variable metrics \cite{Chouzenoux2016}, which makes use of specific preconditioning matrices. For more general nonconvex optimization problems, it is interesting to find the KL exponents of such objective functions in order to find their local convergence rates. Extra theoretical and empirical work in these directions is expected in the future. 

\subsubsection*{Acknowledgments}
The authors would like to thank Jian-Feng Cai and Jinshan Zeng for useful discussions and comments. 
\bibliographystyle{abbrvnat}
\bibliography{nips_2017_ref}
\addcontentsline{toc}{section}{References}

\appendix
\section*{Appendix}
\addcontentsline{toc}{section}{Appendix}
In this Appendix, we provide more details not mentioned in the main text due to space constraint. 

\begin{algorithm}[H]
			\caption{APGnc with adaptive momentum (\APGncp)}
			\label{apgncp}
			\SetKwInOut{Input}{Input}
			\Input{$\by_1=\bx_1=\bx_0, \beta, t\in(0, 1), \eta\le\nicefrac{1}{L}$. }
			\For{$k=1, 2, \ldots$}{
				$\bx_{k}=\prox_{\eta g}(\by_k-\eta\nabla f(\by_k))$,\;
				$\bv_k=\bx_k+\beta(\bx_k-\bx_{k-1})$,\;
				\If{$F(\bx_k)\le F(\bv_k)$}{$\by_{k+1}=\bx_k$, $\beta\leftarrow t\beta$,}
				\ElseIf{$F(\bx_k)\ge F(\bv_k)$}{$\by_{k+1}=\bv_k$, $\beta\leftarrow \min\{\beta/t, 1\}$.}
			}
\end{algorithm}

\begin{algorithm}[H]
	\caption{Block prox-linear (BPL)}
	\label{bpl}
	\SetKwInOut{Input}{Input}
	\Input{$\bx^{(1)}=\bx^{(0)}=\bx^{(-1)}$. }
	\For{$k=1, 2, \ldots$}{
		Pick $b_k\in\{1, \ldots, s\}$ in a deterministic or random manner,\;
		Set $\alpha_k>0$ and $\omega_k\ge 0$,\;
		$\begin{dcases*}
		\bx_i^{(k)} = \bx_i^{(k-1)} & if $i\ne b_k$,\\
		\bx_i^{(k)}\in\prox_{\alpha_kg_i}\left(\widehat{\bx}_i^{(k)} -\alpha_k\nabla_{\bx_i} f\left(\bx_{\ne i}^{(k-1)}, \widehat{\bx}_i^{(k)}  \right)  \right) & if $i=b_k$,\\
		\text{where }\widehat{\bx}_i^{(k)} =\bx_i^{(k-1)}+\omega_k(\bx_i^{(k-1)}-\bx_i^{(k-2)}).
		\end{dcases*}$		
	}
	
\end{algorithm}

\begin{notation}
	$\bx_{< i}$ is the shorthand notation for $(\bx_1, \ldots, \bx_{i-1})$, $\bx_{> i}$ is the short-hand notation for $(\bx_{i+1}, \ldots, \bx_{s})$, $\bx_{\ne i}$ means $(\bx_{< i}, \bx_{> i})$. So $f(\bx_{\ne i}, \widehat{\bx}_i)$ means $f(\bx_{< i}, \widehat{\bx}_i, \bx_{> i})$. 
\end{notation}

\begin{definition}[Proximity operator \cite{Combettes2011}] 
Let $\lambda$ be a positive parameter. The \emph{proximity operator} $\prox_{\lambda g}:\RR^n\to\RR^n$ is defined through
	\[\prox_{\lambda g} (\bx):=\argmin_{\by\in\RR^n} g(\by)+\dfrac{1}{2\lambda}\|\by-\bx\|^2 . \]
	If $g$ is convex, proper and lower semicontinuous, $\prox_g$ admits a unique solution. If $g$ is nonconvex, then it is generally set-valued. 
\end{definition}
\begin{definition}[Domain]
The \emph{domain} of $F: \RR^n\to\RR
$ is defined by
	\[\dom F:=\{\bx\in\RR^n:F(\bx)<+\infty\}. \]
\end{definition}

\begin{definition}[Subdifferential \cite{Rockafellar1998}] \hfill
	\begin{enumerate}
		\item For a given $\bx\in\dom F$, the \emph{Fr\'{e}chet subdifferential} of $F$ at $\bx$, written $\widehat{\partial}F(\bx)$, is the set of all vectors $\bu\in\RR^n$ which satisfy
		\[\lim\limits_{\by\ne \bx}\inf_{\by\to \bx}\dfrac{F(\by)-F(\bx)-\langle \bu, \by-\bx\rangle}{\|\by-\bx\|}\ge0. \]
		When $\bx\notin\dom F$, we set $\widehat{\partial}F(\bx)=\emptyset$. 
		\item The \emph{limiting-subdifferential}, or simply the \textit{subdifferential}, of $F$ at $\bx\in\dom F$, written $\partial F(\bx)$, is defined through the following closure process
		\[\partial F(\bx):=\{\bv\in\RR^n:\exists \bx^k\to x, F(\bx^k)\to F(\bx), \bv^k\in\widehat{\partial}F(\bx^k)\to \bv\}. \]
	\end{enumerate}
\end{definition}

\begin{definition}[Sublevel set] 
Being given real numbers $\alpha$ and $\beta$ we set 
	\[[\alpha\le F\le \beta]:=\{\bx\in\RR^n: \alpha\le F(\bx)\le \beta\}.\]
	$[\alpha<F<\beta]$ is defined similarly. 
\end{definition}

\begin{definition}[Distance]
	The distance of a point $\bx\in\RR^n$ to a closed set $\Omega\subseteq\RR^n$ is defined as
	\[\mathrm{dist}(\bx, \Omega):=\inf_{\by\in\Omega} \|\by-\bx\|.\]
	If $\Omega=\emptyset$, we have that $\mathrm{dist}(\bx,\Omega)=\infty$ for all $\bx$. 
\end{definition}

\begin{definition}[Kurdyka-\L{}ojasiewicz property and KL function \cite{Bolte2014}] \label{klprop}\hfill
	\begin{enumerate}
		\item The function $F:\RR^n\to\RR\cup\{+\infty\}$ is said to have the \emph{Kurdyka-\L{}ojasiewicz property} at $\bar{\bx}\in\dom\partial F$ if there exist $\eta\in(0, +\infty]$, a neighbourhood $\calB_\rho(\bar{\bx}):=\{\bx:\|\bx-\bar{\bx}\|<\rho\}$ of $\bar{\bx}$ and a continuous concave function $\varphi(t):=ct^{1-\theta}$ for some $c>0$ and $\theta\in[0,1)$ such that for all $\bx\in \calB_\rho(\bar{\bx})\cap[F(\bar{\bx})<F<F(\bar{\bx})+\eta]$, the Kurdyka-\L{}ojasiewicz inequality holds
			\[\varphi'(F(\bx)-F(\bar{\bx}))\,\mathrm{dist}(\bm{0}, \partial F(\bx))\ge1. \] 
		\item Proper lower semincontinuous functions which satisfy the Kurdyka-\L{}ojasiewicz inequality at each point of $\dom\partial F$ are called \textit{KL functions}. 
	\end{enumerate}
\end{definition}

\begin{definition}[KL exponent \cite{Li2017b}] \label{klexp}
	For a proper closed function $F$ satisfying the KL property at $\bar{\bx}\in\dom\partial F$, if the corresponding function $\varphi$ can be chosen as in \Cref{klprop}, the KL inequality can be written as 
	\[\mathrm{dist}(\bm{0}, \partial F(\bx))\ge \bar{c}(F(\bx)-F(\bar{\bx}))^\theta \]
	for some $\bar{c}>0$. We say $F$ has the KL property at $\bar{\bx}$ with an exponent $\theta$. If $F$ is a KL function and has the same exponent $\theta$ at any $\bar{\bx}\in\dom\partial F$, then we can say that $F$ is a KL function with an exponent of $\theta$. 
\end{definition}

\begin{definition}[SCAD penalty \cite{Fan2001}]\label{scad}
The SCAD penalty $r_{\lambda, \gamma}(u)$ is defined as 
\begin{equation*}
 r_{\lambda, \gamma}(u) =\begin{cases}
 \lambda|u|, & \text{if }|u|<\lambda,\\
 \frac{2\gamma\lambda|u|-(u^2+\lambda^2)}{2(\gamma-1)}, & \text{if }\lambda<|u|\le\gamma\lambda, \\
 \frac{\lambda^2(\gamma^2-1)}{2(\gamma-1)}, & \text{if }|u|>\gamma\lambda. 
 \end{cases}
\end{equation*}
\end{definition}

\begin{proposition}[Proximal-Gradient Gauss-Southwell rules \cite{Gong2013}]\hfill\em
\begin{enumerate}
\item In coordinate descent methods, the GS-$s$ rule chooses the coordinate with the  most negative directional derivative, given by
\[i_k=\argmax_i\left\lbrace \min_{s\in\partial g_i}|\nabla_i f(\bx^{(k)})+s|\right\rbrace. \]
We generalize it to the block coordinate scenario which has the form
\[i_k=\argmax_i\left\lbrace \min_{\bs\in\partial g_i}\|\nabla_{\bx_i} f(\bx^{(k)})+\bs\|\right\rbrace. \]
\item The GS-$r$ rule selects the coordinate which maximizes the length of the step 
\[i_k=\argmax_i\left\lbrace\left| x_i^{(k)}-\prox_{g_i/L}\left(x_i^{(k)}-\dfrac{1}{L}\nabla_i f(\bx^{(k)}) \right) \right|  \right\rbrace, \]
which is generalized to the block coordinate version
\[i_k=\argmax_i\left\lbrace\left\| \bx_i^{(k)}-\prox_{g_i/L}\left(\bx_i^{(k)}-\dfrac{1}{L}\nabla_{\bx_i} f(\bx^{(k)}) \right) \right\|  \right\rbrace.\]
\item The GS-$q$ rule maximizes the progress assuming a quadratic upper bound on $f$ 
\[i_k=\argmax_i\left\lbrace \min_d\left[ f(\bx^{(k)})+\nabla_i f(\bx^{(k)})d+\dfrac{L}{2}d^2+g_i(\bx_i^{(k)}+d)-g_i(\bx_i^{(k)})\right] \right\rbrace. \]
We do not apply this to a block coordinate scenario. 
\end{enumerate}
\end{proposition}

\begin{proposition}[Proximity operators of regularizers in \S3]\em
If $\bx, \bu\in\RR^n$, then
\begin{enumerate}
\item for $R(\bx)=\lambda\|\bx\|_1$ \cite{Bach2012}, 
\[\prox_{\lambda\|\cdot\|_1}(\bu)= \left(\prox_{\lambda|\cdot|}(\bu_i) \right)_{1\le i \le n} =\left(\sign(\bu_i)(|\bu_i|-\lambda)_+\right)_{1\le i \le n},\]
where $(u)_+=\max\{0, u\}$ and the above proximity operator is often referred to as \emph{soft-thresholding}. 

\item for the $\ell_1/\ell_2$-norm $R(\bx)=\lambda\|\bx\|_2$ \cite{Bach2012},
\[\prox_{\lambda\|\cdot\|_2}(\bu)=\left(1-\dfrac{\lambda}{\|\bu\|_2} \right)_+\bu.  \]
Further, if $\calG$ is a partition of $\{1, \ldots, n\}$, for the $\ell_1/\ell_2$-norm $R:\bx\mapsto\lambda\sum_{g\in\calG} \|\bx_g\|_2$ we have 
\[\left( \prox_{R}(\bu)\right)_g= \left(1-\dfrac{\lambda}{\|\bu_g\|_2} \right)_+\bu_g, \]
which is often referred to as \emph{group-soft-thresholding}. The problem being solved is called \textit{group Lasso} if it is used with a least square loss in \Cref{reg_ls}. 

\item for $r_{\lambda, \gamma}(x)$ defined in \Cref{scad} \cite{Gong2013}, let $h_1(x)=\frac{1}{2}(x-u)^2+r_{\lambda, \gamma}(x)$ and 
\begin{align*}
u_1 &:= \sign(u)\min\{\lambda, (|u|-\lambda)_+\},\\
u_2 &:= \sign(u)\min\left\lbrace  \gamma\lambda, \max\left\lbrace \lambda, \frac{|u|(\gamma-1)-\gamma\lambda}{\gamma-2}\right\rbrace \right\rbrace ,\\
u_3 &:= \sign(u)\max\{\gamma\lambda, |u|\}. 
\end{align*}
Then, we have
\[\prox_{r_{\lambda, \gamma}}(u)=\begin{cases}
u_1 & \text{if }u_1=\argmin_y h_1(y), \\
u_2 & \text{if }u_2=\argmin_y h_1(y), \\
u_3 & \text{if }u_3=\argmin_y h_1(y).
\end{cases}\]
Thus, we further have 
\[\prox_{\sum_{i=1}^nr_{\lambda, \gamma}}(\bu)=\left(\prox_{r_{\lambda, \gamma}}(\bu_i) \right)_{1\le i\le n}. \]
\item for $R(\bx)=\lambda\sum_{i=1}^n\min\{|x_i|,\theta\}$ \cite{Gong2013}, let $h_2(x)=\frac{1}{2}(x-v)^2+\lambda\min\{|x|,\theta\}$ and 
\begin{align*}
v_1 &:= \sign(v)\max\{\theta, |v|\},\\
v_2 &:= \sign(v)\min\{\theta, (|v|-\lambda)_+\}.
\end{align*}
Then, we have 
\[\prox_{\lambda\min\{|\cdot|,\theta\}}(v)=\begin{cases}
v_1 & \text{if }h_2(v_1)\le h_2(v_2),\\
v_2 & \text{otherwise.}
\end{cases}\]
Thus, we further have 
\[\prox_{R}(\bv)=\left(\prox_{\lambda\min\{|\cdot|,\theta\}}(\bv_i) \right)_{1\le i\le n}. \]
\end{enumerate}
\end{proposition}

\end{document}